\documentclass[leqno, 11pt]{amsart}
\usepackage{amsmath, amssymb, latexsym}
  \newtheorem{definition}{Definition}[section]
  
  \newtheorem{lemma}[definition]{Lemma}
  \newtheorem{proposition}[definition]{Proposition}

  \theoremstyle{remark}
  
  \newtheorem{remark}[definition]{Remark}

  \newtheorem*{acknowledgements}{Acknowledgements}

  \numberwithin{equation}{section}

\newcommand{\op}[1]{\operatorname{#1}}

\newcommand{\C}{\ensuremath{\mathbb{C}}} 
\newcommand{\N}{\ensuremath{\mathbb{N}}} 
\newcommand{\R}{\ensuremath{\mathbb{R}}} 
\newcommand{\Z}{\ensuremath{\mathbb{Z}}} 
\newcommand{\bC}{\ensuremath{\overline{\C}_{-}}}
\newcommand{\bCt}{\ensuremath{\overline{\C}_{-,\tau}}}

\newcommand{\cS}{\ensuremath{\mathcal{S}}}

\newcommand{\psido}{$\Psi$DO} 
\newcommand{\pvdo}{\Psi_{\op{v}}}
\newcommand{\psidos}{$\Psi$DO's}

\newcommand{\hotimes}{\hat\otimes}

\newcommand{\vo}{\textup{v}}
\newcommand{\voa}{\textup{v}, \textup{a}}
\newcommand{\vob}{\scriptscriptstyle{\|}\vo}

\begin{document}
    \title{On the Asymptotic completeness\\ of the Volterra calculus}
\author{Rapha\"el Ponge \medskip \\
With an Appendix by H.~Mikayelyan and R.~Ponge}

\address{\noindent Department of Mathematics, Ohio State University, Columbus, OH 43210,~USA.}
\email{ponge@math.ohio-state.edu}

\address{\emph{Address of Hayk Mikayelyan:}
Mathematics Institute, University of Leipzig ~~\\ \indent Augustusplatz 10/11, D-04109 Leipzig, Germany.}
\email{hayk@math.uni-leipzig.de}

\thanks{The author was partially supported by the European RT Network \emph{Geometric Analysis} HPCRN-CT-1999-00118.}



\begin{abstract}
    The Volterra calculus is a simple and powerful pseudodifferential tool for inverting parabolic equations and it has 
    also found many applications in geometric analysis. On the other hand, an important property in the theory of pseudodifferential 
    operators is the asymptotic completeness, which allows us to construct parametrices modulo smoothing 
    operators. In this paper we present new and fairly elementary proofs the asymptotic completeness of the 
    Volterra calculus. 
\end{abstract}

\maketitle 

\section*{Introduction}
This paper deals with the asymptotic completeness of the Volterra calculus.  
Recall that the 
latter was invented in the early 70's by Piriou~\cite{Pi:COPDTV} and Greiner~\cite{Gr:AEHE} and  consists in a 
modification of the classical \psido\ calculus in order to take into account two classical properties occurring in the 
context of parabolic equations: the Volterra property and the anisotropy with respect to the time variable 
(cf.~Section~\ref{sec.Volterra}). As a consequence  the Volterra calculus proved to be a powerful tool for inverting parabolic equations 
(see Piriou~(\cite{Pi:COPDTV}, \cite{Pi:PLGODPDB})) and for deriving small heat kernel asymptotics for elliptic 
operators~(see Greiner~\cite{Gr:AEHE}).   

Subsequently, the Volterra calculus has  been extended to several other settings. In~\cite{BGS:HECRM} Beals-Greiner-Stanton 
produced a version of the Volterra calculus for the 
hypoelliptic calculus on Heisenberg manifolds~(\cite{BG:CHM}, \cite{Ta:NCMA}) and used it to derive the small 
time heat kernel asymptotics for the Kohn Laplacian on CR manifolds. Also, Melrose~\cite{Me:APSIT} 
fit the Volterra calculus into the framework of his $b$-calculus on manifolds with boundary 
and used it to invert the heat equation with the purpose of producing a heat kernel proof of 
the Atiyah-Patodi-Singer index theorem~\cite{APS:SARG}. 

More recently,  
Buchholz-Schulze~\cite{BS:}, Krainer~(\cite{Kr:PhD}, \cite{Kr:}) and    
Krainer-Schulze~\cite{KS:}  extended the Volterra calculus to the setting of the cone  
calculus of Schulze~(\cite{Sc1:}, \cite{Sc2:}) in order to solve general parabolic problems 
on manifold with conical singularities and to deal with large time asymptotics of solutions
 to parabolic 
 problems on manifolds with boundary 
 (by looking at the infinite time cylinder as a manifold with a conical 
 singularity at time $t=\infty$; see~\cite{Kr:PhD}, \cite{KS:}). 
Furthermore, Mitrea~\cite{Mi:} used a version of the Volterra calculus for 
studying parabolic equations with 
Dirichlet boundary conditions on Lipschitz domains and Mikayelyan~\cite{Mik:PhD} dealt with parabolic 
problems on manifolds with edges via an extension of the Volterra calculus to the setting of Schulze's edge 
calculus~(\cite{Sc1:}, \cite{Sc2:}).

On the other hand, in~\cite{Po:NSPLIFSA} the approach to the heat kernel asymptotics of 
Greiner~\cite{Gr:AEHE} was combined with the rescaling of Getzler~\cite{Ge:SPLASIT} to produce  
a new short proof of the local index formula of Atiyah-Singer~\cite{AS:IEO3}. The upshot is that this proof  
is as simple as Getzler short proof in~\cite{Ge:SPLASIT} but, unlike the latter, it allows us to similarly compute the 
Connes-Moscovici cocycle~\cite{CM:LIFNCG} for 
Dirac spectral triples.  Furthermore, the pseudodifferential representation of the heat kernel provided by the 
Volterra calculus in~\cite{Gr:AEHE} also gives an alternative to the construction by Seeley~\cite{Se:CPEO} of 
pseudodifferential complex powers of (hypo)elliptic differential operators (\emph{cf.}~\cite{Po:PhD}, \cite{Po:MAMS1}; 
see also~\cite{MSV:FAI}).      

While most of the usual properties of the classical \psido\ calculus hold \emph{verbatim} in the setting of 
the Volterra calculus, a more delicate issue is to check asymptotic completeness. 
This property allows us to construct parametrices for parabolic operators, but its standard proof 
cannot be carried through in the setting of the Volterra calculus. Indeed, at the  level of 
symbols the Volterra property corresponds to analyticity with respect to the 
time covariable (see~Section~\ref{sec.Volterra}), but this property is not preserved by the cut-off arguments 
of the proof. 

Since we cannot cut off Volterra symbols, 
Piriou~\cite[pp.~82--88]{Pi:COPDTV}  proved the 
asymptotic completeness of the Volterra calculus by cutting off distribution 
kernels instead, which at this level does not harm the Volterra property,
and by checking that under the 
Fourier transform we get an actual asymptotic expansion of symbols (see also~\cite{Me:APSIT}). Recently, 
Krainer~\cite[pp.~62--73]{Kr:}
obtained  a proof by making use of the kernel cut-off operator of Schulze~(\cite{Sc1:}, \cite{Sc2:}) and 
Mikayelyan~\cite{Mik:ASOVVS} produced another proof by combining translations in the time covariable 
with an induction process\footnote{Despite that in~\cite[p.~79]{Mik:ASOVVS} 
the induction hypothesis is not stated properly and there is 
a typo on line 14 the argument in the proof is correct.}.

In this paper, we present somewhat simpler approaches. First, we 
show that we actually get  a 
 Volterra \psido\ by adding a suitable smoothing operator to the \psido\ 
 provided by the standard proof of the  asymptotic completeness of classical symbols (see 
 Proposition~\ref{prop:asymptotic-completeness1}).    
 
 Second, we deal with the asymptotic completeness of analytic Volterra symbols (see 
 Proposition~\ref{prop:asymptotic-completeness2} and Proposition~\ref{prop:asymptotic-completeness3}). This was the setting 
 under consideration in~~\cite{Kr:} and~\cite{Mik:ASOVVS}, because this asymptotic completeness implies that 
 of the Volterra calculus (see Section~\ref{sec.proof2}).  Here our approach is 
 inspired by the version of the Borel lemma for analytic functions on an angular sector (e.g.~\cite[p.~63]{AG:}).
 
 This paper is organized as follows. In Section~\ref{sec.Volterra} we briefly review the main facts concerning the 
 Volterra calculus. In Section~\ref{sec.proof1} we present our first approach. In Section~\ref{sec.proof2} 
 we carry out our proofs of the asymptotic completeness of analytic symbols. Finally, in the appendix,  
 written with Hayk Mikayelyan, we give alternative proofs of the asymptotic completeness of these analytic Volterra 
 symbols  by combining our approach with the use of translations in the time covariable from~\cite{Mik:ASOVVS}. 
 In particular we remove the induction process used in that paper.
 

\begin{acknowledgements} 
  { \small I would like to thank Thomas Krainer and Hayk Mikayelyan for useful discussions about the Volterra calculus, as well as   
  the whole PDE group of Profs.~Elmar Schrohe and  Bert-Wolfgang Schulze at Potsdam University for its warm hospitality. }
\end{acknowledgements}

\section{Overview of the Volterra calculus} 
\label{sec.Volterra}
Throughout this paper $U$ is an open subset of $\R^{n}$ and $w$ denotes  an even integer~$\geq 2$. Also,  
we let $\C_{-}$ denote 
the half-space $\C_{-}=\{\Im \tau<0\}\subset \C$ with closure $\bC=\{\Im \tau \leq 0\}$. 

As alluded to in the introduction the Volterra calculus is a pseudodifferential calculus on $U\times \R$ which 
aims to take into account: \smallskip 

(i) The anisotropy of parabolic problems on $U\times \R$, i.e.~their homogeneity with respect to the 
dilations of $\R^{n}\times \bC$ given by 
\begin{equation}
    \lambda.(\xi,\tau)=(\lambda\xi,\lambda^{w}\tau) , \qquad \lambda\in 
        \R\setminus 0.
\end{equation}

(ii) The Volterra property, that is the fact for a continuous operator $Q$ from 
$C^{\infty}_{c}(U_{x}\times\R_{t})$ to 
$C^{\infty}(U_{x}\times\R_{t})$ to have a distribution kernel of the form 
$k_{Q}(x,t; y,s)=K_{Q}(x,y,t-s)$, where $K_{Q}(x,y,t)$ vanishes in the region $U\times U\times\{t<0\}$. 

\begin{definition}\label{def:homogeneous-Volterra-symbols}
    $S_{\vo, m}(U\times \R^{n+1})$, $m\in \Z$, consists of smooth functions $q_{m}(x,\xi,\tau)$ on 
     $U_{x}\times (\R^{n+1}_{(\xi,\tau)}\setminus0)$ such that $q_{m}(x,\xi,\tau)$ can be extended to a smooth 
     function on $U_{x}\times[(\R^{n}_{\xi}\times \bCt)\setminus 0]$ in such way to be analytic with respect to 
     $\tau\in \C_{-}$ 
     and to be homogeneous of degree $m$, i.e.~$q_{m}(x,\lambda 
     \xi,\lambda^{2}\tau)=\lambda^{m}q_{m}(x,\xi,\tau)$ for any $\lambda \in \R\setminus 0$. 
\end{definition}

In fact, Definition~\ref{def:homogeneous-Volterra-symbols} is  
intimately related to the Volterra property, since we have:  

\begin{lemma}[{\cite[Prop.~1.9]{BGS:HECRM}}]
    Any symbol $q(x,\xi,\tau)\in S_{\vo, m}(U\times\R^{n+1})$ can be extended into a unique distribution 
    $g(x,\xi,\tau) \in 
    C^{\infty}(U)\hotimes \cS'(\R^{n+1})$ in such way to be homogeneous with respect to the covariables 
    $(\xi,\tau)$ and such that $\check{q}(x,y,t):=\mathcal{F}^{-1}_{(\xi,\tau)\rightarrow (y,t)}[g](x,y,t)$ vanishes for $t<0$. 
\end{lemma}

Next, we introduce the pseudo-norm on $\R^{n}\times \bC$ given by 
\begin{equation}
    \|\xi,\tau\|=(|\xi|^{w}+|\tau|)^{1/w}, \qquad  (\xi,\tau)\in \R^{n}\times \bC.
\end{equation}
This pseudo-norm is homogeneous since $\|(\lambda \xi,\lambda^{w}\tau)\|=|\lambda|^{w}\|\xi,\tau\|$ for any 
$\lambda\in \R\setminus 0$. Also, for $ (\xi,\tau)\in \R^{n}\times \bC$ we have
\begin{equation}
2^{-1/w}(1+|\xi|+|\tau|)^{1/w}\leq 1+\|\xi,\tau\|\leq 1+|\xi|+|\tau|.
     \label{eq:estimates-pseudo-norm}
\end{equation}

\begin{definition}%
$S_{\vo}^m(U\times\R^{n+1})$, $m\in\Z$,  consists of smooth functions $q(x,\xi,\tau)$ on 
    $U_{x}\times\R^{n+1}_{(\xi,\tau)}$ which have an asymptotic expansion  $q \sim \sum_{j\geq 0} q_{m-j}$, 
    where $q_{m-j}\in 
    S_{\vo, m-j}(U\times \R^{n+1})$ and  $\sim$ means that, for any integer $N\geq 0$ and 
    for any compact $K\subset U$,  we have 
    \begin{equation}
                      |\partial^{\alpha}_{x}\partial^{\beta}_{\xi} \partial^k_{\tau}(q-\sum_{j<N} 
            q_{m-j})(x,\xi,\tau) | 
               \leq C_{NK\alpha\beta k} \|\xi,\tau\|^{m-|\beta|-wk-N},
                          \label{eq:volterra.asymptotic-symbols}
            \end{equation}
            for $x\in K$ and for $(\xi,\tau) \in \R^{n+1}$ such that $\|\xi,\tau\|\geq 1$.
\end{definition}

\begin{remark}\label{rem:Volterra.type0,1/p}
    It follows from~(\ref{eq:estimates-pseudo-norm}) and~(\ref{eq:volterra.asymptotic-symbols}) that 
    $S_{\vo}^{m}(U\times\R^{n+1})$ is contained in the H\"ormander's class $S^{m'}_{0,\frac{1}{w}}((U_{x}\times 
    \R_{t})\times \R^{n+1}_{(\xi,\tau)})$ with $m'=m$ if $m\geq 0$ and $m'=\frac{m}{w}$ otherwise.  In fact, using H\"ormander's 
    Lemma~(e.g.~\cite[Thm.~2.9]{Ho:}, \cite[Prop.~3.6]{Sh:POST}) one can even show that the asymptotic expansion in the sense 
    of~(\ref{eq:volterra.asymptotic-symbols}) coincides with that for standard symbols. 
\end{remark}

\begin{definition}\label{def:Volterra.PsiDO's}
    $\pvdo^m(U\times\R)$, $m\in\Z$,  consists of continuous operators 
    $Q$ from $C_{c}^\infty(U\times\R)$ to 
    $C^\infty(U\times\R)$ such that:\smallskip 
   
   (i)  $Q$  has the Volterra property;\smallskip 
   
   (ii) $Q$ is of the form $Q=q(x,D_{x},D_{t})+R$ 
    with $q\in S^m_{\vo}(U\times\R^{n+1})$ and $R$ smoothing operator. 
\end{definition}

As  it follows from Remark~\ref{rem:Volterra.type0,1/p} the class $\pvdo^{*}(U\times\R)$ is contained in 
the class of \psidos\ of type 
$(0,\frac{1}{w})$ on $U\times \R$. 
Therefore, once the asymptotic completeness is checked, all the standard properties of classical \psidos\ 
hold \emph{verbatim} for Volterra \psidos: symbolic calculus, existence of parametrices for parabolic \psidos\ 
(i.e.~ those with an invertible principal symbol), invariance by diffeomorphisms which don't act on the time 
variable. 
In particular, the Volterra calculus makes sense on $M\times \R$ for any smooth manifold $M$. 

On the other hand, the Volterra calculus has two important applications:\smallskip 

- \emph{Inversion of parabolic operators} (Piriou~(\cite{Pi:COPDTV}, \cite{Pi:PLGODPDB})). Any parabolic 
differential operator on $U\times \R$, not only admits a parametrix, but has actually an inverse in the 
Volterra calculus. This makes use of the well known 
fact that if $R$ is a smoothing operator which is properly supported and has the Volterra property,  then the Levi series 
$\sum_{j \geq 1} R^{j}$ is convergent in the Fr\'echet space of smoothing operators. This result has been extended to 
several  other settings (see \cite{BGS:HECRM}, \cite{Me:APSIT}, \cite{BS:}, \cite{Kr:PhD}, \cite{Kr:}, \cite{KS:}, 
\cite{Mik:PhD}, \cite{Mi:}).\smallskip 

- \emph{Heat kernel asymptotics}~(Greiner~\cite{Gr:AEHE}). Let $P$ be differential operator of order $w$ on a compact 
Riemannian manifold $M$ and assume that the principal symbol of $P$ is positive definite. Then we can relate 
the heat kernel $k_{t}(x,y)$ of $P$ to the the distribution kernel of $(P+\partial_{t})^{-1}$ so that, as the latter 
is a Volterra \psido, we can derive the asymptotics for $k_{t}(x,x)$ as $t\rightarrow 0^{+}$ in terms of the 
symbol of $(P+\partial_{t})^{-1}$. As alluded to in the introduction this approach to the heat kernel asymptotics has been extended 
to the setting 
of the hypoelliptic calculus on Heisenberg manifolds~(see~\cite{BGS:HECRM}) and has been used for proving the local index formula of 
Atiyah-Singer~(cf.~\cite{Po:NSPLIFSA}) and for constructing complex powers of (hypo)elliptic 
operators~(cf.~\cite{Po:PhD},  \cite{Po:MAMS1}; see also~\cite{MSV:FAI}).

\section{Asymptotic completeness of the Volterra calculus}
\label{sec.proof1}
Here we give our first proof of the asymptotic completeness of the Volterra calculus. More precisely, we 
shall prove: 

\begin{proposition}\label{prop:asymptotic-completeness1}
  Given $q_{m-j}\in S_{\vo, m-j}(U\times\R^{n+1})$, $j=0,1,\ldots$, there always exists 
   $Q\in \Psi^{m}_{\vo}(U\times \R)$ with symbol $q\sim \sum_{j\geq 0}q_{m-j}$. 
\end{proposition}
\begin{proof}
      For $\epsilon >0$ and $(\xi,\tau) \in \R^{n} \times \R$ 
      let $c_{\epsilon}(\xi,\tau) =1- \phi(\epsilon \|\xi,\tau\|)$, 
      where $\phi(u)\in C^{\infty}_{c}([0, \infty))$ is such that $\phi(u)=1$ near $u=0$. 
      Then similar arguments as those in the standard proof of the asymptotic completeness of symbols 
      (e.g.~\cite[Thm.~2.7]{Ho:}, \cite[Prop.~2.5]{Sh:POST})  show that for any $\epsilon \geq 1$ and for any compact 
      $K\subset U$  we have 
      \begin{equation}
          |\partial_{x}^{\alpha}\partial_{\xi}^{\beta}  \partial_{\tau}^{k}
          [c_{\epsilon}(\xi,\tau)q_{m-j}(x,\xi,\tau)]| \leq C_{jK\alpha\beta k} 
           \epsilon (1+\|\xi,\tau\|)^{m+1-j-|\beta|-wk},
          \label{eq:proof1.estimates}
      \end{equation}
       for $(x,\xi,\tau)\in K \times \R^{n}\times \R$ and where the constant $C_{jK\alpha\beta k}$ does depend 
       on $\epsilon$. 
       
       Next, given $(K_{j})_{j \geq 0}$ an increasing compact exhaustion of $U$ the 
       estimates~(\ref{eq:proof1.estimates}) allows us to find numbers $\epsilon_{j}\geq 1$, 
       $j=0,1,\ldots$, such that
      \begin{equation}
          |\partial_{x}^{\alpha}\partial_{\xi}^{\beta} 
          \partial_{\tau}^{k}[c_{\epsilon_{j}}(\xi,\tau) q_{m-j}(x,\xi,\tau)]| \leq 2^{-j}
          (1+\|\xi,\tau\|)^{m+1-j-|\beta|-wk},  
          \label{eq:proof1.estimates-2j}
      \end{equation}
       for $l+|\alpha|+\beta|+k<j$ and $(x,\xi,\tau)\in K_{l}\in  \R^{n}\times \R$. Therefore, the series 
       $\sum_{j\geq 0} c_{\epsilon_{j}}  q_{m-j}$ converges in 
       $C^{\infty}(U\times \R^{n+1})$ to some function $q$. Moreover, the estimates~(\ref{eq:proof1.estimates-2j}) also imply that 
       $q\sim \sum_{j\geq 0} q_{m-j}$. Thus, $q$ is in $S_{\vo}^{m}(U\times \R^{n+1})$. 
     
     Nevertheless, the operator $q(x,D_{x}, D_{t})$ needs not have the 
     Volterra property, since the cut-off functions $c_{\epsilon_{j}}(\xi,\tau)$ kill the 
     analyticity of $q_{m-j}(x,\xi, \tau)$ with respect to $\tau$. Thus, 
     we  need to construct a smoothing operator $R$ such that $q(x,D_{x}, D_{t})+R$ 
    has the Volterra property. 
    
    First, as the Fourier transform relates the decay at infinity to the behavior at the origin 
    of the Fourier transform, 
     the estimates~(\ref{eq:volterra.asymptotic-symbols}) imply that for any integer $N$ the 
     distribution $\check{q}(x,y,t) -\sum_{j\leq J} \check{q}_{m-j}(x,y,t)$ is in 
     $C^{N}(U_{x}\times\R^{n}_{y}\times\R_{t})$ as soon as $J$ is large enough. As $\check{q}_{m-j}(x,y,t)$ vanishes for 
     $t<0$ it follows that for every  integer $l\geq 0$ the limit $\lim_{t \rightarrow 0^-} \partial_{t}^l \check{q}(.,.,t)$ exists in 
     $C^N(U\times\R^n)$ for any $N \geq l$, hence exists in $C^\infty(U\times\R^n)$. 
     
     Now, using a version of the Borel lemma with coefficients in the Fr\'echet space $C^{\infty}(U\times\R^{n})$ 
     we can construct a smooth function 
     $R(x,y,t)$ on $U\times\R^{n}\times\R$ such that for any integer $l\in \N$ we have 
     $\partial_{t}^l R(.,.,0)= 
              \lim_{t\rightarrow 0^-}\partial_{t}^l 
              \check{q}(.,.,t)$ in $C^{\infty}(U\times\R^{n})$.
      Then on $U\times \R^{n} \times \R$ we define 
      \begin{equation}
                R_{1}(x,y,t)=(1-\chi(t))(\check{q}(x,y,t)-R(x,y,t)), 
                \label{eq:}
       \end{equation}
     where $\chi(t)$ denotes the characteristic function of the interval $[0,\infty)$. In fact, $R_{1}(x,y,t)$ is a 
     smooth function on $U\times \R^{n}\times \R$. Indeed, $R_{1}(x,y,t)$ is obviously smooth for $t\neq 0$ 
     and, as 
     $\partial_{t}^l R_{1}(.,.,t)=0$ for $t>0$ and as we have 
     $\lim_{t\rightarrow 0^-} \partial_{t}^l R_{1}(x,y,t)=0$ in $C^{\infty}(U\times \R^{n})$, we see 
     that $R_{1}(x,y,t)$ is also smooth near $t=0$. 
     
     Finally, let $Q:C^\infty_{c}(U\times\R^{n+1})\rightarrow C^\infty(U\times\R^{n+1})$ be 
     the operator with distribution  kernel 
     \begin{multline}
         K_{Q}(x,y,t-s)=\chi(t-s)(\check{q}(x,x-y,t-s)-R(x,x-y,t-s)),\\
         = \check{q}(x,x-y,t-s)-R(x,x-y,t-s) - R_{1}(x,x-y,t-s).  
     \end{multline}
      Then $Q$ has the Volterra property and differs from $q(x,D_{x},D_{t})$ by a smoothing operator, so
     is a Volterra \psido\ with symbol $q\sim \sum_{j\geq 0} q_{m-j}$.
 \end{proof}

\section{Asymptotic completeness of analytic Volterra symbols}
\label{sec.proof2}
Using a different approach, partly inspired by the proof of the Borel lemma for analytic functions on an 
angular sector (see~\cite[p.~63]{AG:}), we will now prove the asymptotic completeness of the analytic 
Volterra symbols below. 

\begin{definition}%
 $S^{m}_{\voa}(U\times\R^{n+1})$, $m\in \Z$, consists of smooth functions $q(x,\xi,\tau)$ on 
    $U_{x}\times\R^{n+1}_{(\xi,\tau)}$ such that: \smallskip 
     
     (i) $q(x,\xi,\tau)$ extends to a smooth function on $U_{x}\times \R^{n}_{\xi}\times \bCt$ in such 
     way to be analytic with respect to $\tau\in \C_{-}$;\smallskip 
     
     (ii) We have $q\sim_{\textup{a}} \sum_{j\geq 0} q_{m-j}$, $q_{m-j} \in S_{\vo, m-j}(U\times\R^{n+1})$, in the sense 
     that, for any integer $N \geq 0$ and for any compact $K \subset U$, we have 
         \begin{equation}
                      |\partial^{\alpha}_{x}\partial^{\beta}_{\xi} \partial^k_{\tau}(q-\sum_{j<N} 
            q_{m-j})(x,\xi,\tau) | 
               \leq C_{NK\alpha\beta k} \|\xi,\tau\|^{m-|\beta|-wk-N},
                           \label{eq:volterra.asymptotic-symbols-analytic}
            \end{equation}
            for $x\in K$ and for $(\xi,\tau) \in \R^{n}\times \bC$ such that $\|\xi,\tau\|\geq 1$.
\end{definition}

In fact, by the Paley-Wiener theorem if $q(x,\xi,\tau) \in S^{m}_{\voa}(U\times\R^{n+1})$ then $\check{q}(x,y,t)=0$ 
for $t<0$. Thus, the operator $q(x,D_{x},D_{t})$ is already a Volterra \psido\ since its distribution kernel  is 
$\check{q}(x,x-y,s-t)$. Thus, the asymptotic completeness of analytic Volterra symbols implies the asymptotic 
completeness of the Volterra calculus. 
 
Next, consider the homogeneous symbol $\rho(\xi,\tau)\in S_{\vo, -1}(\R^{n+1})$ given by 
\begin{equation}
    \rho(\xi,\tau)=(|\xi|^{p}+i\tau)^{-1/w}, \qquad (\xi,\tau)\in (\R^{n}\times \bC)\setminus 0,
\end{equation}
where in order to define the $w$'th root we use the continuous determination of the argument on $\C\setminus[0,-\infty)$
with values in 
$(-\pi,\pi)$, so that $\rho(\xi,\tau)$ takes values in $\Omega=\{z\in \C\setminus0; \ |\arg z|\leq 
\frac{\pi}{2w}\}$. Moreover, as $\rho(\xi,\tau)$ never vanishes on $(\R^{n}\times \bC)\setminus 0$ and is homogeneous of degree 
$-1$ there exists $C_{\rho}>0$ such that for $(\xi,\tau)\in (\R^{n}\times \bC)\setminus 0$ we have
\begin{equation}
    C^{-1}_{\rho}\|\xi,\tau\|^{-1}\leq \rho(x,\xi)\leq C_{\rho} \|\xi,\tau\|^{-1}. 
     \label{eq:estimates-rho}
\end{equation}

Now,  for any integer $N$ we have $z^{N}e^{-z}\rightarrow 0$ as $z\in \Omega$ goes to infinity. Therefore, for 
any $\epsilon>0$ we define a smooth function on $\R^{n}\times \bC$ by letting 
\begin{equation}
 a_{\epsilon}(0,0)=0 \quad \text{and} \quad   a_{\epsilon}(\xi,\tau)= e^{-\epsilon \rho(\xi,\tau)} \ 
 \text{for $(\xi,\tau)\neq 0$}. 
\end{equation}
Notice that $a_{\epsilon}(\xi,\tau)$ is analytic with respect to $\tau \in \C_{-}$. In fact, we have: 

\begin{lemma}\label{lem:a-epsilon}
    1) $a_{\epsilon}$ is in $S^{0}_{\voa}(\R^{n+1})$ and we have 
    $a_{\epsilon}\sim_{\textup{a}} \sum_{j\geq 0} \frac{\epsilon^{j}}{j!} \rho^{j}$.\smallskip 
    
    2) For any $\epsilon \geq 1$ and any integer $N\geq 0$ we have 
    \begin{gather}
        |\partial_{\xi}^{\beta}\partial_{\tau}^{k} a_{\epsilon}(\xi,\tau)| \leq C_{\beta k} \epsilon^{-1} 
        \|\xi,\tau\|^{1-|\beta|-wj}, \qquad \|\xi,\tau\| \geq 1, 
        \label{eq:estimates.a-epsilon1}\\
         |\partial_{\xi}^{\beta}\partial_{\tau}^{k} a_{\epsilon}(\xi,\tau)| \leq C_{N\beta k} \epsilon^{-1} 
         \|\xi,\tau\|^{N}, \qquad \|\xi,\tau\| \leq 1,
         \label{eq:estimates.a-epsilon2}
    \end{gather}
    where the constants $C_{\beta j}$ and $C_{N\beta j}$ are independent of $\epsilon$. 
\end{lemma}
\begin{proof}
    First,  if $\|\xi,\tau\|\geq 1$ then by~(\ref{eq:estimates-rho}) we have $\rho(\xi,\tau)\leq 
    C_{\rho}$, and so we get:    
    \begin{equation}
        |a_{\epsilon}(x,\tau)-\sum_{j<J} \frac{\epsilon^{j}}{j!}\rho(\xi,\tau)^{j}|\leq 
        |\rho(x,\xi)|^{J}\sum _{j\geq J}\frac{\epsilon^{j}}{j!} C_{\rho}^{j-J}\leq C_{\epsilon J}  
        \|\xi,\tau\|^{-J}.
         \label{eq:asymptotic.a-epsilon}
    \end{equation}
     
    On the other hand, an easy induction shows that for any multi-order $\beta$ and any integer $j$ the 
    function $\partial_{\xi}^{\beta}\partial_{\tau}^{k} a_{\epsilon}(\xi,\tau)$ is a linear combination of 
    terms of the form $ \epsilon^{l} \eta_{\beta kl}(\xi,\tau) e^{-\epsilon \rho(\xi,\tau)}$,
     where $l$ is an integer $\leq j$ and $\eta_{\beta kl}(\xi,\tau)$ is  homogeneous of degree $-(|\beta|+wk)-l$ 
    and does not depend on $\epsilon$.  In particular, as $\epsilon\geq 1$ and as for 
    any $N\geq 0$ the function $z^{N}e^{-z}$ is bounded on $\Omega$, we get 
    \begin{multline}
         \|\xi,\tau\|^{-N} \epsilon^{l}  |\eta_{\beta kl}(\xi,\tau)e^{-\epsilon \rho(\xi,\tau)}|=\\ 
         \epsilon^{-(N+1)}.\|\xi,\tau\|^{-N}  |\eta_{\beta kl}(\xi,\tau)\rho(\xi,\tau)^{-(N+l+1)}|.|(\epsilon 
         \rho(\xi,\tau))^{N+l+1}e^{-\epsilon \rho(\xi,\tau)}|,\\
         \leq C_{\beta 
         klN}\epsilon^{-1}\|\xi,\tau\|^{1+|\beta|+wk},
    \end{multline}
    where the constant $C_{\beta klN}$ does not depend on $\epsilon$. 
    Then by setting $N=0$ we obtain~(\ref{eq:estimates.a-epsilon1}) and 
    by taking $N$ large enough we get~(\ref{eq:estimates.a-epsilon2}). 
     
    Finally, thanks to the H\"ormander Lemma~(\cite[Thm.~2.9]{Ho:}, \cite[Prop.~3.6]{Sh:POST}) the 
    estimates~(\ref{eq:estimates.a-epsilon1})--(\ref{eq:asymptotic.a-epsilon}) are enough to show that 
    $a_{\epsilon}\sim_{\textup{a}} \sum_{j\geq 0}  \frac{\epsilon^{j}}{j!} \rho^{j}$. In particular,  
    the symbol $a_{\epsilon}$ belong to $S^{0}_{\vo}(\R^{n+1})$. 
\end{proof}

\begin{proposition}\label{prop:asymptotic-completeness2}
    For $j=0,1,2, \ldots$ let $q_{m-j}\in S_{\vo, m-j}(U\times \R^{n+1})$. Then there exists 
    $q\in S^{m}_{\vo, a}(U\times \R^{n+1})$ such that $q\sim_{\textup{a}} \sum_{j \geq 0} q_{m-j}$. 
    In particular, the operator 
    $q(x,D_{x},D_{t})$ is a Volterra \psido\ with symbol $q\sim \sum_{j \geq 0} q_{m-j}$. 
\end{proposition}
\begin{proof}
    We seek for numbers  $\epsilon_{j}\geq 1$ and symbols $r_{m-j} \in S_{\vo, m-j}(U\times 
    \R^{n+1})$, $j=0,1,\ldots$, such that: \smallskip 
    
    (i) The series $\sum_{j \geq 0} a_{\epsilon_{j}} (\xi,\tau)r_{m-j}(x,\xi,\tau)$ converges in 
    $C^{\infty}(U\times \R^{n}\times \bC)$ to some function $q(x,\xi,\tau)$ which is analytic with respect 
    to $\tau \in \bC$; \smallskip 
    
    (ii) We have $q \sim_{\textup{a}} \sum_{j\geq 0}a_{\epsilon_{j}}r_{m-j}$. \smallskip 
     
   \noindent  Notice that by Lemma~\ref{lem:a-epsilon} 
the function $a_{\epsilon_{j}}(\xi,\tau)r_{m-j}(x,\xi,\tau)$ is smooth on $U\times 
\R^{n}\times \bC$ and analytic with respect to $\tau \in \C_{-}$, so that (i) makes sense. Also, 
Lemma~\ref{lem:a-epsilon} implies that 
$a_{\epsilon}r_{m-j} \sim_{\textup{a}} \sum_{k\geq 
    0}\frac{\epsilon_{j}^{k}}{k!}\rho^{k}r_{m-j}$. Therefore, if (ii) holds then we 
    obtain 
\begin{equation}
    q \sim_{\textup{a}} \sum_{j \geq 0} a_{\epsilon}r_{m-j} \sim_{\textup{a}} \sum_{j,l\geq 0} 
    \frac{\epsilon_{j}^{l}}{l!}\rho^{l}r_{m-j}.
     \label{eq:proof2.asymptotic-a-rmj}
\end{equation}
Thus, we would have $q \sim_{\textup{a}} \sum_{j\geq 0} q_{m-j}$ if, and only if, for $j=0,1,\ldots$ we have 
\begin{equation}
    q_{m-j}=r_{m-j}+\epsilon_{j-1}\rho r_{m-j+1}+\ldots+ 
    \frac{\epsilon_{0}^{j}}{j!}\rho^{j}r_{m}, \qquad j\geq 0.
    \label{eq:proof2.asymptotic-a-rmj2}
\end{equation}
By an easy induction these equalities allow us to uniquely determine $r_{m-j}$ in terms of $q_{m},\ldots,q_{m-j}$ and 
$\epsilon_{0},\ldots,\epsilon_{j-1}$ only, so that $r_{m-j}$ does not depend on $\epsilon_{l}$ for $l\geq 0$.  
Therefore, using~(\ref{eq:estimates.a-epsilon1}) and (\ref{eq:estimates.a-epsilon2}) we see that for any compact 
$K\subset U$    we have 
\begin{equation}
     |\partial_{x}^{\alpha}\partial_{\xi}^{\beta}\partial_{\tau}^{k} [a_{\epsilon_{j}}
     r_{m-j}](x,\xi,\tau)| \leq C_{K\alpha\beta kj}
     \epsilon_{j}^{-1}(1+\|\xi,\tau\|)^{m-j+1-|\beta|-wk}, 
     \label{eq:estimates-a-epsilon-rm-j}
\end{equation}
for $x\in K$ and for $(\xi,\tau)\times \R^{n}\times \bC$.

Now, let $(K_{j})_{j \geq 0}$ be an increasing exhaustion of $U$ by compact subsets. Then thanks 
to~(\ref{eq:estimates-a-epsilon-rm-j})  
we can choose the sequence $(\epsilon_{j})_{j\geq 0}$ in such way that we have 
\begin{equation}
     |\partial_{x}^{\alpha}\partial_{\xi}^{\beta}\partial_{\tau}^{k} [a_{\epsilon_{j}} 
     r_{m-j}](x,\xi,\tau)]| \leq 2^{-j}(1+ \|\xi,\tau\|)^{m-j+1-|\beta|-wk},
     \label{eq:estimates-a-epsilon-rm-j-2j}
\end{equation}
for $l+|\beta|+k\leq j$ and $(x,\xi,\tau)\in K_{l}\times \R^{n}\times \bC$. 

It follows from~(\ref{eq:estimates-a-epsilon-rm-j-2j}) that the series $\sum_{j\geq 0} a_{\epsilon_{j}}(\xi,\tau)r_{m-j}(x,\xi,\tau)$ 
converges in $C^{\infty}(U\times  \R^{n}\times\bC)$ to some function $q(x,\xi,\tau)$. This function is 
furthermore is analytic with respect 
to $\tau \in \C_{-}$ since each term $a_{\epsilon_{j}}(\xi,\tau)r_{m-j}(x,\xi,\tau)$ in the series is. 

On the other hand, the 
estimates~(\ref{eq:estimates-a-epsilon-rm-j-2j}) also imply that 
$q\sim_{\textup{a}} \sum_{j\geq 0} a_{\epsilon_{j}}r_{m-j}$, which in view of~(\ref{eq:proof2.asymptotic-a-rmj}) 
and~(\ref{eq:proof2.asymptotic-a-rmj2}) yields $q\sim_{\textup{a}}\sum_{j\geq 0} 
q_{m-j}$. In particular, the function $q$ belongs to $S^{m}_{\vo, a}(U\times \R^{n+1})$. 
\end{proof}

This approach also allows us to deal with the asymptotic completeness of non-polyhomogeneous analytic Volterra 
symbols. These symbols can be defined as follows. 
\begin{definition}
 $S^{m}_{\vob}(U\times\R^{n+1})$, $m\in \R$, consists of smooth functions $q(x,\xi,\tau)$ 
    on  $U\times\R^{n+1}_{(\xi,\tau)}$ which can be extended to a smooth function on $U\times \R^{n}\times 
    \bC$ in such way that: \smallskip 
    
    (i) $q(x,\xi,\tau)$ is analytic with respect to $\tau\in \C_{-}$; \smallskip 
   
   (ii) For any compact $K \subset U$ we have 
    \begin{equation}
        |\partial_{x}^{\alpha}\partial_{\xi}^{\beta}\partial_{\tau}^{k} q(x,\xi,\tau)| \leq C_{K\alpha\beta 
        k} (1+\|\xi,\tau\|)^{m-|\beta|-w|k},
    \end{equation}  
    for $x\in K$ and $(\xi,\tau)\in \R^{n}\times \bC$. 
\end{definition}

\begin{remark}
 Any symbol in $S_{\voa}(U\times\R^{n+1})$ is contained in 
 $S^{m}_{\vob}(U\times\R^{n+1})$. Furthermore, if $q\in 
 S^{m_{j}}_{\vob}(U\times\R^{n+1})$ where $m_{j}\downarrow -\infty$ as $m_{j} 
 \rightarrow \infty$ then we have $q \sim_{\textup{a}} \sum_{j\geq 0}$ if, and only if, for any integer $N\geq 0$ the 
 symbol  $q-\sum_{j\leq J} q_{j}$ is $S^{-N}_{\vob}(U\times\R^{n+1})$ as soon as 
 $J$ is large enough.
\end{remark}

\begin{proposition}\label{prop:asymptotic-completeness3}
  Let $q_{j}\in S^{m_{j}}_{\vo \scriptscriptstyle{\|}}(U\times\R^{n+1})$, $j=0,1,\ldots$, where 
    \mbox{$m_{j}\downarrow-\infty$} as $j\rightarrow \infty$. Then there exists 
    $q\in S^{m_{0}}_{\vob}(U\times\R^{n+1})$ such that 
    $q\sim_{\textup{a}} \sum_{j\geq 0} q_{j}$. 
\end{proposition}
\begin{proof}
    First, we can always assume $m_{j}-1\leq m_{j+1}$ for any $j\geq 0$, possibly 
    by replacing the sequence 
    $(q_{j})_{j \geq 0}$ 
    by the sequence $(q_{j,l})$, which is indexed by couples $(j,l)\in \N^{2}$ such that $0\leq j\leq 
    m_{j}-m_{j+1}$ and is given by 
    \begin{equation}
        q_{j,0}=q_{j} \quad \textup{and} \quad  q_{j,l}=0 \ \textup{for $1\leq l\leq m_{j}-m_{j+1}$}.
    \end{equation}
    This has the effect to insert finitely many zero terms of order $\geq m_{j+1}$ into the sequence 
    $(q_{j})_{j \geq 0}$, so does not affect the class of symbols that are asymptotic to $\sum_{j \geq 0} 
    q_{j}$.  
    
    Bearing this assumption in mind we now seek for numbers $\epsilon_{j} \geq 1$ and symbols $r_{j} \in 
    S^{m_{j}}$, $j=0,1,..$, such that: \smallskip 
    
    (i) The series $\sum_{j\geq 0} a_{\epsilon_{j}}(\xi,\tau)r_{j}(x,\xi,\tau)$ 
converges in $C^{\infty}(U\times  \R^{n}\times\bC)$ to a function $q(x,\xi,\tau)$ which is analytic with 
respect to $\tau \in \C_{-}$; \smallskip 

(ii) We have $q \sim_{\textup{a}} \sum_{j\geq 0} a_{\epsilon_{j}}r_{j}$. \smallskip 

\noindent  As in~(\ref{eq:proof2.asymptotic-a-rmj2}) the 
condition~(ii) would imply that $q \sim_{\textup{a}} \sum_{j\geq 0} q_{j}$ if we choose 
the symbols $r_{j}$ in such way that for $j=0,1,\ldots$ we have 
\begin{equation}
    q_{j}=\!\!\!\sum_{m_{j+1} <m_{k}-l\leq m_{j}} \!\!\! \frac{\epsilon_{k}^{l}}{l!}\rho^{l}r_{k}
    = r_{j} +\!\!\!\sum_{\substack{m_{j+1} <m_{k}-l\leq m_{j},\\ k>j}}\!\!\! \frac{\epsilon_{k}^{l}}{l!}\rho^{l}r_{k},
     \label{eq:proof2.asymptotic-a-rmj3}
\end{equation}
where the second equality holds because $m_{j}-1\leq m_{j+1}$.  
This uniquely determines $r_{j}$ in terms of $q_{0},\ldots,q_{j}$ and $\epsilon_{0}, 
\ldots,\epsilon_{j-1}$ only. Therefore, along the same lines as that of the proof of 
Proposition~\ref{prop:asymptotic-completeness2} we can find numbers $\epsilon_{j} \geq 1$ 
 such that (i) and (ii) hold. Then thanks to~(\ref{eq:proof2.asymptotic-a-rmj3}) we have $q\sim_{\textup{a}} \sum_{j\geq 0} q_{j}$. 
\end{proof}
\begin{remark}%
   For some authors (Buchholz-Schulze~\cite{BS:}, Krainer~(\cite{Kr:PhD},\cite{Kr:}), 
   Krainer-Schulze~\cite{KS:}, Mikayelyan~\cite{Mik:PhD}) 
   the Volterra \psidos\ are defined as those coming from analytic Volterra symbols 
   only. Since Proposition~\ref{prop:asymptotic-completeness2} implies that  any 
   Volterra \psido\ in the sense of Definition~\ref{def:Volterra.PsiDO's} 
   coincides up to a  smoothing  operator  with the quantification of an analytic Volterra symbol, 
   it follows that the two possible definitions are actually equivalent.
\end{remark}

\appendix 
\section*{Appendix by H.~Mikayelyan and R.~Ponge}
\setcounter{section}{1}
\setcounter{equation}{0}
\setcounter{definition}{0}
In this appendix we present alternative proofs of the asymptotic 
completeness of the analytic Volterra symbols by combining the use of  translations in the time covariable from~\cite{Mik:ASOVVS} 
with some of the ideas from Section~\ref{sec.proof2}.  
In particular we remove the induction process used in~\cite{Mik:ASOVVS}.

In the sequel given a symbol $q$ on $U\times\R^{n}\times \bC$ 
for any $T>0$ we let 
\begin{equation}
    q^{(T)}(x,\xi,\tau)=q(x,\xi,\tau-iT). 
\end{equation}
We shall first deal with non-polyhomogeneous symbols, which is the setting under consideration in~\cite{Mik:ASOVVS}. 

\begin{lemma}[Krainer~(\cite{Kr:PhD}, \cite{Kr:})]\label{lem:Appendix.Krainer}
 If $q \in S^{m}_{\vob}(U\times\R^{n+1})$, $m\in \R$, then 
 the symbol $q^{(T)}$ is in $S_{\vob}^{m}(U\times\R^{n+1})$  and we have 
 $q^{(T)} \sim_{\textup{a}} \sum_{l\geq 0} \frac{(-iT)^{l}}{l!}\partial_{\tau}^{l} q$. 
\end{lemma}
\begin{proof}
Since $T>0$  we have $|\tau|\leq |\tau-iT|\leq |\tau|+T$ for any $\tau\in \bC$. Therefore, 
    for any $(\xi,\tau) \in \R^{n}\times \bC$ we have 
    \begin{equation}
        1+\|\xi,\tau\| \leq 1+\|\xi,\tau-iT\| \leq (1+T)^{1/w}(1+\|\xi,\tau\|). 
         \label{eq:Appendix.estimates}
    \end{equation}
If we combine these inequalities with a Taylor formula about $\tau=0$ then for any compact $K\subset U$ we get   
    \begin{multline}
             |\partial_{x}^{\alpha}\partial_{\xi}^{\beta}\partial_{\tau}^{k}  
             [q^{(T)}-\sum_{l=0}^{N} \frac{(-iT)^{l}}{l!} \partial_{\tau}^{l}q](x,\xi,\tau)|   \\
           \leq   C_{TNK\alpha\beta k} \int_{0}^{1} (1+\|\xi,\tau- isT\|)^{m-|\beta|-kw-N-1}ds, \\ 
          \leq C_{TNK\alpha\beta k} (1+\|\xi,\tau\|)^{m-|\beta|-kw-N-1},
         \label{eq:Appendix.Taylor-formula}
    \end{multline}
    for $x\in K$ and $(\xi,\tau)\in \R^{n}\times \bC$. Hence $q^{(T)} \sim_{\textup{a}} \sum_{l\geq 0} 
    \frac{(-iT)^{l}}{l!}\partial_{\tau}^{l} q$. In particular, the function $q$ belongs to 
    $S_{\vob}^{m}(U\times\R^{n}\times \bC)$. 
\end{proof}

\begin{lemma}[Mikayelyan~\cite{Mik:ASOVVS}]\label{lem:Appendix.Mikayelyan}
   Let $q\in S_{\vob}^{m}(U\times \R^{n+1})$ with $m\leq -1$. Then for any compact $K 
    \subset U$ we have 
    \begin{equation}
         |\partial_{x}^{\alpha}\partial_{\xi}^{\beta}\partial_{\tau}^{k} q^{(T)}(x,\xi,\tau)| \leq C_{K\alpha\beta 
        k} (1+T)^{-1/w}(1+\|\xi,\tau\|)^{m+1-|\beta|-w|k},
    \end{equation}
     for $x\in K$ and $(\xi,\tau)\in \R^{n}\times \bC$ and where the constant $C_{K\alpha\beta 
        k}$ does not depend on $T$. 
\end{lemma}
\begin{proof}
  Since $m \leq -1$ and since for $\tau \in \bC$ we have $|\tau-iT| \geq \sup(T, |\tau|)$, we see that 
  for  $(\xi,\tau)\in \R^{n}\times \bC$ we get: 
 \begin{gather}
     (1+\|\xi,\tau-iT\|)^{-1} \leq (1+T)^{-1/w}, \label{Appendix.estimates-Tw}\\ 
      (1+\|\xi,\tau-iT\|)^{m+1-|\beta|-kw} \leq (1+\|\xi,\tau\|)^{m+1-|\beta|-kw}.
 \end{gather}
 Therefore, for any compact $K \subset U$ we have 
  \begin{equation}
    \begin{split}
       |\partial_{x}^{\alpha}\partial_{\xi}^{\beta}\partial_{\tau}^{k} q^{(T)}(x,\xi,\tau)| & \leq C_{K\alpha\beta k} 
       (1+\|\xi,\tau-iT\|)^{-1+m+1-|\beta|-wk},\\
        &\leq C_{K\alpha\beta k} (1+T)^{-1/w}(1+\|\xi,\tau\|)^{m+1-|\beta|-wk},
    \end{split}
  \end{equation}
  for $x\in K$ and $(\xi,\tau)\in \R^{n}\times \bC$ and where the constants $C_{K\alpha\beta 
        k}$ do not depend on $T$. 
\end{proof}

We can now give a second proof of Proposition~\ref{prop:asymptotic-completeness3}. 

\begin{proof}[Second proof of Proposition~\ref{prop:asymptotic-completeness3}] Here we let 
    $q_{j}\in S^{m_{j}}_{\vo \scriptscriptstyle{\|}}(U\times\R^{n+1})$, $j=0,1,\ldots$, with 
    \mbox{$m_{j}\downarrow-\infty$} as $j\rightarrow \infty$ and we look 
    for $q\in S^{m_{0}}_{\vob}(U\times\R^{n+1})$ such that 
    $q\sim_{\textup{a}} \sum_{j\geq 0} q_{j}$. 
    
    First, as in the first proof of Proposition~\ref{prop:asymptotic-completeness3} we can assume $m_{j}-w\leq m_{j+1}$, 
    possibly by replacing 
    the sequence $(q_{j})_{j \geq 0}$ 
    by the sequence $(q_{j,l})$ which is indexed by the couples $(j,l)\in \N^{2}$ such that $0\leq j\leq 
    w^{-1}(m_{j}-m_{j+1})$ and is given by  $q_{j,0}=q_{j}$ if $l=0$ and  $q_{j,l}=0$ if $1\leq l\leq m_{j}-m_{j+1}$.
    
    Bearing this assumption in mind we now seek  for numbers $T_{j}>0$ and  
    symbols $r_{j}\in S_{\vob}^{m_{j}}(U\times \R^{n+1})$, $j=0,1,2,..$, such that: \smallskip 
    
    (i) The series $\sum_{j=0}^{\infty} r_{j}^{(T_{j})}(x,\xi, \tau)$ converges  in $C^{\infty}(U\times 
    \R^{n}\times \bC)$ to some symbol $q(x,\xi, \tau)$ which is analytic with respect to $\tau \in \C_{-}$;\smallskip 
    
    (ii) We have $q\sim_{\textup{a}} \sum_{j \geq 0} r_{j}^{(T_{j})}$. \smallskip 
    
    \noindent Assuming that (i) and (ii) hold, using Lemma~\ref{lem:Appendix.Krainer} we get 
    \begin{equation}
        q\sim_{\textup{a}} \sum_{j\geq 0} r_{j}^{(T_{j})}\sim_{\textup{a}} 
            \sum_{j,l\geq 0} \frac{(-iT_{j})^{l}}{l!}\partial_{\tau}^{l} r_{j}.
     \end{equation}
    Therefore, we would have $q\sim_{\textup{a}} \sum_{j\geq 0} q_{j}$ if  for $j=0,1,\ldots$ we have 
    \begin{equation}
        q_{j}= \!\!\! \sum_{m_{j+1}< m_{j'}-lw \leq m_{j}} \!\!\!  \frac{(-iT_{j'})^{l}}{l!}\partial_{\tau}^{l} r_{j'} = r_{j}+ 
         \!\!\! \sum_{\substack{m_{j+1}< m_{j'}-lw \leq m_{j},\\ j'<j}} \!\!\!  \frac{(-iT_{j'})^{l}}{l!}\partial_{\tau}^{l} r_{j'},
         \label{eq:Appendix.qj-rj}
    \end{equation}
    where the second equality holds because $m_{j}-w \leq m_{j+1}$. 
    
    Now, we set $T_{j}=1$ for all indices such that $m_{j}> -1$. Since the equalities~(\ref{eq:Appendix.qj-rj}) uniquely determine  
    $r_{j}$ in terms of $q_{0}, \ldots, q_{j}$ and $T_{0}, \ldots, T_{j-1}$ only, when $m_{j} \geq -1$  
    Lemma~\ref{lem:Appendix.Mikayelyan} implies that for any compact $K\subset U$ we have 
    \begin{equation}
         |\partial_{x}^{\alpha}\partial_{\xi}^{\beta}\partial_{\tau}^{k} r^{(T_{j})}_{j}(x,\xi,\tau)| \leq C_{K\alpha\beta 
        kj} (1+T_{j})^{-1/w}(1+\|\xi,\tau\|)^{m_{j}+1-|\beta|-wk},
    \end{equation}
    for $x\in K$ and $(\xi, \tau)\in \R^{n}\times \bC$. Then by similar arguments as those of the proof of 
    Proposition~\ref{prop:asymptotic-completeness2}   
    we can construct a sequence of positive numbers $(T_{j})_{m_{j} \leq -1}$ converging fast enough to $\infty$ 
    such that the condition (i) and (ii) above hold. Then~(\ref{eq:Appendix.qj-rj}) implies that $q 
    \sim_{\textup{a}} \sum_{j\geq 0}q_{j}$. 
 \end{proof}

 Next, we deal with the polyhomogeneous case. 
 
 \begin{lemma}\label{lem:Appendix.Krainer-Mikayelyan2}
    Let $q \in S_{\vo,m}(U\times \R^{n+1})$. Then: \smallskip 
    
    1) $q^{(T)}$ belongs to $S^{m}_{\vo, a}(U\times \R^{n+1})$ and we have 
    $q^{(T)} \sim_{\textup{a}} \sum_{l\geq 0} \frac{(-iT)^{l}}{l!}\partial_{\tau}^{l} q$. \smallskip 
    
    2) If $m \leq -1$ and $T\geq 1$ then for any compact $K \subset U$ we have 
    \begin{equation}
         |\partial_{x}^{\alpha}\partial_{\xi}^{\beta}\partial_{\tau}^{k} q^{(T)}(x,\xi,\tau)| \leq C_{K\alpha\beta 
        k} T^{-1/w}(1+\|\xi,\tau\|)^{m+1-|\beta|-wk},
    \end{equation}
     for $(x,\xi,\tau)\in K\times\R^{n}\times \bC$ and where the constant $C_{K\alpha\beta 
        k}$ does not depend on $T$. 
 \end{lemma}
 \begin{proof}
     First, as $T$ is positive $\bC -iT$ is contained in $\C_{-}$. Since $q(x,\xi,\tau)$ is smooth and 
     analytic with respect to $\tau$ on $U \times\R^{n}\times \C_{-}$, it follows that 
     $q^{(T)}(x,\xi,\tau)$  is smooth on $U\times \R^{n}\times \bC$  and is analytic with respect to $ \tau \in \C_{-}$. 
     
     Next, as in~(\ref{eq:Appendix.Taylor-formula}) by combining a Taylor formula with~(\ref{eq:Appendix.estimates}) we get   
     \begin{multline}
        |\partial_{x}^{\alpha}\partial_{\xi}^{\beta}\partial_{\tau}^{k}( q^{(T)}-\sum_{l\leq N} 
        \frac{(-iT)^{l}}{l!}\partial_{\tau}^{l} q)(x,\xi,\tau)|   \\
          \leq   C_{TNK\alpha\beta k} \int_{0}^{1} \|\xi,\tau- isT\|^{m-|\beta|-kw-N-1}ds, \\ \leq 
            C_{TNK\alpha\beta k} \|\xi,\tau\|^{m-|\beta|-kw-N-1},
     \end{multline}
     for $x \in K$ and for $(\xi,\tau) \in \R^{n}\times \bC$ such that $\|\xi,\tau\|\geq 1$. Thus $q^{(T)}$ is asymptotic to 
     $\sum_{l \geq 0 }\frac{(-iT)^{l}}{l!}\partial_{\tau}^{l} q$ in the sense of~(\ref{eq:volterra.asymptotic-symbols-analytic}).   
     Hence $q^{(T)}$ belongs to 
     $S^{m}_{\vo,a}(U\times \R^{n+1})$. 
     
     On the other hand, as in~(\ref{Appendix.estimates-Tw}) since $|\tau -iT|\geq T$ we have 
    \begin{equation}
         \|\xi,\tau-iT\|^{-1} \leq T^{-1/w}.
    \end{equation}
     Moreover, if  $T\geq 1$ then for $(\xi,\tau) \in \R^{n}\times \bC$ we also get 
     \begin{equation}
         \|\xi,\tau-iT\| \geq [|\xi|^{w}+(|\Re \tau|^{2}+(|\Im \tau|+T)^{2})^{1/2}]^{1/w} \geq C_{w} 
         (1+\|\xi,\tau\|).
     \end{equation}
Therefore, for any $T \geq 1$ and for any compact $K \subset      U$ we have 
\begin{equation}
    \begin{split}
               |\partial_{x}^{\alpha}\partial_{\xi}^{\beta}\partial_{\tau}^{k} q^{(T)}(x,\xi,\tau)| & 
               \leq C_{K\alpha\beta k} \|\xi,\tau-iT\|^{-1+m+1-|\beta|-wk}, \\
               & \leq C_{K\alpha\beta k} T^{-1/w}(1+\|\xi,\tau\|)^{m+1-|\beta|-wk},
    \end{split}
\end{equation}
     for $x\in K$ and $(\xi,\tau)\in \R^{n}\times \bC$ and where the constant $C_{K\alpha\beta k}$ does not depend 
     on $T$. 
\end{proof}
 
\begin{proof}[Second proof of Proposition~\ref{prop:asymptotic-completeness2}] 
    For $j=0,1,\ldots$ let $q\in S_{\vo, m-j}(U\times\R^{n+1})$. Then, provided that we make use of 
    Lemma~\ref{lem:Appendix.Krainer-Mikayelyan2} instead of Lemma~\ref{lem:Appendix.Krainer} and 
    Lemma~\ref{lem:Appendix.Mikayelyan}, similar arguments as those of the second proof of 
    Proposition~\ref{prop:asymptotic-completeness3} show that we can find numbers $T_{j}\geq 1$ and symbols 
    $r_{m-j}\in S_{\vo, m-j}(U\times\R^{n+1})$, $j=0,1,...$, such that: \smallskip 
    
    (i) The series $\sum_{j\geq 0} r_{m-j}^{(T_{j})}$ converges  in $C^{\infty}(U\times 
    \R^{n}\times \bC)$ to some symbol $q(x,\xi, \tau)$ which is analytic with respect to $\tau \in \C_{-}$;\smallskip 
    
    (ii) We have $q\sim_{\textup{a}} \sum_{j \geq 0} r_{m-j}^{(T_{j})}\sim_{\textup{a}} \sum_{j\geq 0} q_{m-j}$. \smallskip 
    
   \noindent  Hence the proposition. 
\end{proof}

\end{document}